\documentclass[11pt]{article}
\usepackage[top=2.5cm,bottom=2.5cm,left=2.5cm,right=2.5cm]{geometry}
\usepackage{mathrsfs}
\usepackage{pstricks}
\usepackage{amsmath,amssymb,graphicx,amsthm,tabularx,array}
\usepackage{hyperref}
\usepackage{makecell}

\theoremstyle{plain}
\newtheorem{theorem}{Theorem}[section]

\newtheorem{coro}[theorem]{Corollary}
\newtheorem{lem}[theorem]{Lemma}

\newtheorem{problem}{Problem}[section]

\theoremstyle{definition}

\newtheorem{other}{}

\title{{A degree-stability method for the $p$-spectral Tur\'an-type  problem for disjoint color-critical graphs}\thanks{
Ya-Lei Jin is partly supported by the National Natural Science Foundation of China (Nos. 12371349 and 12471331),
Ming-Zhu Chen is partly supported by the National Natural Science Foundation of China (No. 12461063),
and Peng-Li Zhang is partly supported by the National Natural Science Foundation of China (No. 12501493).
E-mail addresses: zhengj@yangtzeu.edu.cn (J. Zheng),  yaleijin@shnu.edu.cn (Y.-L. Jin), mzchen@hainanu.edu.cn (M.-Z. Chen), zpengli0626@163.com ($^\dag$P.-L. Zhang, corresponding author).}}

\author{Jian Zheng$^a$, Ya-Lei Jin$^b$, Ming-Zhu Chen$^c$,  Peng-Li Zhang$^d\dag$\\
{\small  $^a$ School of  Information and Mathematics, Yangtze University, Jingzhou 434023,  P.R. China} \\
{\small $^b$  Department of Mathematics, Shanghai Normal University, Shanghai 200234, P.R. China}\\
{\small  $^c$ School of Mathematics and Statistics, Key Laboratory of Engineering Modeling and Statistical } \\
{\small  Computation of Hainan Province,
Hainan University, Haikou 570228, P.R. China}\\
{\small $^d$ School of Statistics and Data Science, Shanghai University of International} \\
{\small Business and Economics, Shanghai 201620, P.R. China}
}

\date{}
\begin{document}
\maketitle
\begin{abstract}
Let $F_1,\ldots,F_t$ be pairwise vertex-disjoint color-critical graphs with $\chi(F_i)=r+1$ for each $i\in[t]$.
For $r\ge 2$ and sufficiently large $n$, Simonovits [Discrete Math.\ 7 (1974) 349--376] proved that the graph $K_{t-1}\vee T_{n-t+1,r}$
uniquely maximizes the number of
edges among $n$-vertex graphs containing no copy of $\cup_{i=1}^tF_i$.
In 2024, Lei and Li [Electron.\ J.\ Combin.\ 31 (2024), P1.25]
established the corresponding
result for the adjacency spectral radius. We prove that, for every fixed $p>1$ and all sufficiently large $n$, the
same graph uniquely maximizes the $p$-spectral radius in this class.
 The proof reduces
the problem to forbidding $t$ disjoint copies of a single color-critical graph and applies degree stability after
deleting $t-1$ such copies.
\end{abstract}



 \section{Introduction}
 A graph $G$ is defined by  a vertex set $V(G)$ and an edge set $E(G)$,
where every edge is a two-element subset of  $V(G)$.
We denote by $e(G):=|E(G)|$  the number of edges of $G$.
For positive integers $a_1,\ldots,a_r$, let $K_r(a_1,\ldots,a_r)$ stand for the complete $r$-partite graph
whose partite sets have sizes  $a_1,\ldots,a_r$.
The \emph{Tur\'{a}n graph} $T_{n,r}$ is the $n$-vertex complete $r$-partite graph
with parts as equal in size as possible.

A \emph{proper vertex coloring} of a graph $F$  is a vertex coloring in which adjacent vertices receive different colors,
and the \emph{chromatic number} $\chi(F)$ is the least number of colors required by such a coloring.
A graph $F$ is  said to be  \emph{color-critical}  if deleting some edge  $e\in E(F)$ strictly reduces the chromatic number, i.e., $\chi(F-e)<\chi(F)$, where $F-e$ denotes the graph obtained from $F$ by  removing $e$.

Given a graph $F$, a graph $G$ is called \emph{$F$-free} if it  contains no copy of $F$ as a subgraph.
The \emph{Tur\'an number} $\mathrm{ex}(n,F)$ is the maximum possible number of edges in an
$F$-free graph on $n$ vertices, and  $\mathrm{EX}(n,F)$ denotes  the family of all $n$-vertex $F$-free graphs attaining this maximum. Tur\'an's classical theorem \cite{Turan41} asserts that $\mathrm{EX}(n,K_{r+1})=\{T_{n,r}\}$.
The  \emph{Tur\'an density} of $F$ is  defined as
 $${\pi}(F):=\lim\limits_{n\to \infty}\frac{ex(n,F)}{\binom{n}{2}},$$
and the celebrated  Erd\H{o}s--Stone--Simonovits theorem \cite{ES1946,ES1966} states that
${\pi}(F)= 1-\frac{1}{\chi(F)-1}$.

For pairwise vertex-disjoint graphs $G_1,\ldots,G_t$, their \emph{union}  $G_1\cup\cdots\cup G_t$
is the graph with vertex set $\cup_{i=1}^{t}V(G_i)$ and edge set $\cup_{i=1}^{t}E(G_i)$.
When $G_1\cong\cdots\cong G_t$, we write $tG_1$ for this union.
The \emph{join} $G_1\vee G_2$ is  obtained from $G_1\cup G_2$
by adding all edges between $V(G_1)$ and $V(G_2)$.

Simonovits \cite{Sim1974} proved the following  Tur\'an result for disjoint color-critical graphs.
\begin{theorem}[\cite{Sim1974}] \label{Sim-col}
Let $F_{1},\ldots, F_{t}$ be pairwise vertex-disjoint color-critical graphs with $\chi(F_{i})=r+1$ for all  $i\in[t]$, where $r\geq2$. Then, for all sufficiently large $n$,  $\mathrm{EX}(n,\cup_{i=1}^{t}F_{i})=\{K_{t-1}\vee T_{n-t+1,r}\}$.
\end{theorem}

Let $G$ be an $n$-vertex graph.  Its \emph{adjacency matrix}  $A(G) = (a_{ij})_{n \times n}$  is defined by $a_{ij} = 1$ if $ij \in E(G)$ and $a_{ij} = 0$ otherwise. The \emph{adjacency spectral radius} $\lambda(G)$ is the largest eigenvalue of $A(G)$. By the Rayleigh-Ritz theorem,
$$\lambda(G)=\max_{\|\mathbf{x}\|_2=1}\mathbf{x}^{T}A\mathbf{x}=\max_{\|\mathbf{x}\|_2=1}2\sum_{uv\in E(G)} x_{u}x_{v}.$$

Correspondingly, we define $\mathrm{ex}_{\lambda}(n,F)$ as the maximum adjacency spectral radius among all $n$-vertex $F$-free graphs, and $\mathrm{EX}_{\lambda}(n,F)$ as the set of $n$-vertex $F$-free graphs attaining this maximum.
The adjacency spectral version of Tur\'an's theorem was established independently by
Nikiforov \cite{Niki2007laa2} and Guiduli \cite{Gui1996}. Later,  Nikiforov \cite{N2009}  proved that $\mathrm{EX}_{\lambda}(n,F) = \{T_{n,r}\}$ for every color-critical graph $F$ with $\chi(F) = r+1$ and sufficiently large $n$.
We refer the reader to \cite{LiLF} for a comprehensive survey on this topic.

Recently, Lei and Li \cite{LL2024} obtained an adjacency spectral analogue of Theorem \ref{Sim-col}.
\begin{theorem}[\cite{LL2024}] \label{LL-col}
Let $F_{1},\ldots, F_{t}$ be pairwise vertex-disjoint color-critical graphs with $\chi(F_{i})=r+1$ for all $i\in[t]$,
where $r\geq2$. Then, for all sufficiently large $n$,
 $\mathrm{EX}_{\lambda}(n,\cup_{i=1}^{t}F_{i})=\{K_{t-1}\vee T_{n-t+1,r}\}$.
\end{theorem}

For a vector $\mathbf{x}=(x_{1},\ldots,x_{n})\in \mathbb{R}^{n}$ and $p \ge 1$, the $\ell_p$-norm of $\mathbf{x}$ is
$\|\mathbf{x}\|_{p}:=(|x_{1}|^{p}+\cdots +|x_{n}|^{p})^{1/p}$.
 For any $p\geq1$, the \emph{$p$-spectral radius} of $G$ is defined by
\begin{equation}\label{def}
\lambda^{(p)}(G)=\max_{\|\mathbf{x}\|_p=1}\mathbf{x}^{T}A\mathbf{x}=\max_{\|\mathbf{x}\|_p=1}2\sum_{uv\in E(G)} x_{u}x_{v}.
\end{equation}
As shown by inequalities $(8)$ and $(9)$ in  \cite{KN2014}, we have
\begin{equation}\label{p-bds}
2e(G)n^{-2/p}\leq \lambda^{(p)}(G)\leq(2e(G))^{1-1/p}.
\end{equation}

\noindent {\bf Remark.} Observe that  $\lambda^{(2)}(G)$ coincides with the adjacency spectral radius, while $\lambda^{(1)}(G)$ is the well-studied Lagrangian of $G$. Moreover, letting $p \to \infty$, $\lambda^{(\infty)}(G)/2$ equals  the number of edges of $G$. Thus $\lambda^{(p)}(G)$ provides a unifying framework for three classical parameters in extremal graph theory.

A vector $\mathbf{x}\in \mathbb{R}^{n}$  is called an \emph{eigenvector} corresponding to $\lambda^{(p)}(G)$
if $\|\mathbf{x}\|_{p}=1$ and $\lambda^{(p)}(G)=\mathbf{x}^{T}A\mathbf{x}$. There is always a nonnegative eigenvector associated with  $\lambda^{(p)}(G)$; we call such a vector a  \emph{principal eigenvector}.
If a principal eigenvector $\mathbf{x}$ is strictly  positive, i.e., $x_{v}>0$ for every $v\in V(G)$,  then it is termed a \emph{Perron--Frobenius eigenvector}.
For $p>1$, a principal eigenvector $\mathbf{x}=(x_{1},\ldots,x_{n})$  satisfies the following eigenequations, obtained via Lagrange's method:
\begin{equation}\label{tzfc}
\lambda^{(p)}(G)x_{k}^{p-1}=\sum_{i\in N_{G}(k)} x_{i},~ ~
1\leq k\leq n.
\end{equation}

We now define  $\mathrm{ex}_{\lambda}^{(p)}(n,F)$ as the maximum $p$-spectral radius among all $n$-vertex $F$-free graphs, and $\mathrm{EX}_{\lambda}^{(p)}(n,F)$ as the set of those graphs attaining this maximum.
The  notion $\lambda^{(p)}(G)$ was first introduced by Keevash, Lenz and Mubayi \cite{KLM2014}, who  established  two key criteria for spectral extremal problems. Applying one of these criteria, they showed that for every  $p>1,$  every color-critical graph  $F$  with $\chi(F)=r+1$ and sufficiently large $n$,  $\mathrm{EX}_{\lambda}^{(p)}(n,F)=\{T_{n,r}\}$. In \cite{KN2014}, Kang and Nikiforov proved the  $p$-spectral version of Tur\'an's theorem and developed a saturation theory for  $\lambda^{(p)}(G)$; a stability theorem for the $p$-spectral radius was later obtained by Li and Peng \cite{LP2022}.
 We remark that extremal problems for $\lambda^{(p)}$ have also been investigated in the hypergraph setting in \cite{KNY2014,NLK2022,N2014A,N2014B,ZLF2024}; nevertheless, the theory for graphs is substantially more mature, and further exploration of extremal questions in the graph case remains worthwhile.


Theorems~\ref{Sim-col}, ~\ref{LL-col} and the above results motivate us to raise the following problem.
\begin{problem}\label{problem}
Let $F_{1},\ldots, F_{t}$ be pairwise vertex-disjoint color-critical graphs with $\chi(F_{i})=r+1$ for all  $i\in[t]$, where $r\geq2$. Is it true that $\mathrm{EX}_{\lambda}^{(p)}(n,\cup_{i=1}^{t}F_{i})=\{K_{t-1}\vee T_{n-t+1,r}\}$ for $p>1$ and   sufficiently large $n$?
\end{problem}

  Let $K^{+}_{r}(a_1, a_2, \dots, a_r)$ be the graph  obtained from a complete $r$-partite graph with parts of sizes $a_1 \geq  2, a_2, \dots, a_r$, by
inserting an additional edge inside the first part. Our main result is as follows.
\begin{theorem}\label{main-thm}
Let $p>1,$ and let $t, a_1, \ldots, a_r$ be positive integers with $a_1 \ge 2$ and $r \ge 2$. Then there exists $n_{0}$ such that for all $n\geq n_{0}$, $\mathrm{EX}_{\lambda}^{(p)}(n,tK^{+}_{r}(a_{1},\ldots,a_{r}))=\{K_{t-1}\vee T_{n-t+1,r}\}$.
\end{theorem}

Let $F_{1},\ldots, F_{t}$ be  any $t$ pairwise vertex-disjoint color-critical graphs with $\chi(F_{i})=r+1$ for all $i\in [t]$.
Clearly, there exist $r$ integers $a_{1},\ldots, a_{r}$ such that $F_{i}$ is a subgraph of $K^{+}_{r}(a_{1},\ldots,a_{r})$ for
every $i\in [t]$.  Since $K_{t-1}\vee T_{n-t+1,r}$ is $\cup_{i=1}^{t}F_{i}$-free,
we obtain the following result from Theorem~\ref{main-thm},
which gives an affirmative answer to Problem~\ref{problem}. 

\begin{coro}
Let $p>1,$ and let $F_{1},\ldots, F_{t}$ be $t\geq 1$ pairwise  vertex-disjoint  color-critical graphs with $\chi(F_{i})=r+1$
for all $i\in[t]$, where $r\geq2$.
Then, for all sufficiently large $n$, $\mathrm{EX}_{\lambda}^{(p)}(n,\cup_{i=1}^{t}F_{i})=\{K_{t-1}\vee T_{n-t+1,r}\}$.
\end{coro}
\paragraph{Proof strategy.}
The degree-stability method is a powerful tool in extremal combinatorics. In this paper, we employ it to study the $p$-spectral Tur\'an-type problem for disjoint color-critical graphs.
Although $K_{t-1}\vee T_{n-t+1,r}$ is $tK^{+}_{r}(a_{1},\ldots,a_{r})$-free and $\delta(K_{t-1}\vee T_{n-t+1,r})\geq (1-\frac{1}{r})n-O(1),$  $\chi(K_{t-1}\vee T_{n-t+1,r})=r+t-1\geq r+1$ for $t\geq 2.$ Hence,
$tK^{+}_{r}(a_{1},\ldots,a_{r})$ is not degree-stable with respect to the family of $r$-partite graphs for $t\geq 2$, which implies that the classical degree-stability method cannot be applied directly. To overcome this difficulty, we adopt a ``delete-then-prove'' strategy. Namely, by removing from the $p$-spectral extremal graph a copy of $(t-1)K^{+}_{r}(a_{1},\ldots,a_{r})$, we obtain a graph that is $K^{+}_{r}(a_{1},\ldots,a_{r})$-free. We then exploit the degree-stability of the color-critical graph
$K^{+}_{r}(a_{1},\ldots,a_{r})$ to restrict the remaining structure to an $r$-partite graph.
Spectral techniques are subsequently employed to precisely control the set of vertices having large degree within their own part, to establish the partite structure of the remaining graph, and ultimately to identify the unique extremal graph
as $K_{t-1}\vee T_{n-t+1,r}$.
In fact, Theorem~\ref{main-thm} holds for  $t\geq 1.$  For $t=1,$ we obtain  Theorem~\ref{main-thm} by applying the  degree-stability method directly.

\paragraph{Organization.}
The remainder of this paper is organized as follows. Section $2$ collects the necessary notation, definitions, and auxiliary results. In Section $3$, we establish a key lemma, Lemma \ref{pro-thm}, which plays an essential role in the proof of Theorem \ref{main-thm}. Section $4$ is devoted to the proof of Theorem \ref{main-thm}.
\section{Preliminaries}

Let $G$ be a graph. For a vertex $v\in V(G)$, the neighborhood of $v$ is
$N_G(v):=\{u\in V(G):uv\in E(G)\}$, and its degree is $d_G(v):=|N_G(v)|$.
The \emph{minimum degree} of $G$ is denoted by $\delta(G)$.
For a subset $U\subseteq V(G)$, the \emph{induced subgraph} $G[U]$ has vertex set $U$
and edge set $\{e\in E(G): e\subseteq U\}$.
We write $N_U(v):=N_G(v)\cap U$ and $d_U(v):=|N_U(v)|$
for the set and the number of neighbors of $v$ in $U$, respectively,
and we set $G-U:=G[V(G)\backslash U]$.

Let $\overline{G}$ be the complement of $G.$
For a set of edges $S\subseteq E(G)$ (resp.\ $S\subseteq E(\overline{G})$), let $G-S$ (resp.\ $G+S$) denote the graph
obtained from $G$ by deleting (resp.\ adding) all edges in $S$.
For $V_1,V_2\subseteq V(G)$, we define
$E_G(V_1,V_2):=\{uv\in E(G): u\in V_1,\ v\in V_2\}$
and $e_G(V_1,V_2):=|E_G(V_1,V_2)|$. When the graph $G$ is clear from the context,
we write $e(V_1,V_2)$ in place of $e_G(V_1,V_2)$.

\begin{lem}[\cite{CFTZ2022}]\label{dx}
Let $V_{1},\ldots, V_{n}$ be $n$ finite sets. Then $|V_{1}\cap\cdots\cap V_{n}|\geq\sum_{i=1}^{n}
|V_{i}|-(n-1)|\cup_{i=1}^{n}V_{i}|$.
\end{lem}

In 2022, Li and Peng \cite[Theorem $3.10$]{LP2022}   presented a  $p$-spectral version of the Erd\H{o}s--Stone--Simonovits theorem.

\begin{theorem}[\cite{LP2022}]\label{p-ess}
Let $F$ be a graph with  $\chi(F)\ge 2$. Then for all $p>1$,
$$\mathrm{ex}_{\lambda}^{(p)}(n,F)=(\pi(F)+o(1))n^{2-2/p}.$$
\end{theorem}

Nikiforov obtained the following lemma via a concise argument; see Claim A in the proof of \cite[Theorem $12$]{N2014A}.
\begin{lem}[\hspace{1sp}\cite{N2014A}]\label{t3}
Let $p>1$, and let $F$ be a graph with $\pi(F)>0$.  Then there exist infinitely many positive integers $n$ satisfying
$$\frac{\mathrm{ex}_{\lambda}^{(p)}(n-1,F)(n-1)^{2/p-1}}{n-2}-\frac{\mathrm{ex}_{\lambda}^{(p)}(n,F)n^{2/p-1}}{n-1}<\frac{1}{n\log n}\cdot\frac{\mathrm{ex}_{\lambda}^{(p)}(n,F)n^{2/p-1}}{n-1}.$$
\end{lem}

\noindent {\bf Fact 1.} Let $p>1$. The function
$$f(x)=\frac{1-2x}{(1-x)^{2/p}}$$ is strictly decreasing on $[0,1)$.

\begin{proof}
A direct differentiation gives
\[f'(x)=\frac{2(1-p+(p-2)x)}{p(1-x)^{2/p+1}}.
\]
If $1<p\leq 2$, then $f'(x)<0$. If $p>2$, then $p-2>0$, and since $x<1
$, we have
\begin{displaymath}
\begin{split}
-(p-1)+(p-2)x<-(p-1)+p-2=-1.
\end{split}
\end{displaymath}
Thus $f'(x)<0$ holds in all cases.
\end{proof}

For a vector $\mathbf{x}=(x_{1},\ldots,x_{n})$, we denote  $x_{\textup{min}}=\min\{x_{1},\ldots,x_{n}\}$.
\begin{lem}\label{t4}
Let $p>1$, and let $F$ be a graph with  $\pi(F)>0$.  There exist infinitely many positive integers $n$ such that, for every $G_n\in\mathrm{EX}_{\lambda}^{(p)}(n,F)$ and every  principal eigenvector
$\mathbf{x}=(x_1,\ldots,x_n)$ corresponding to $\lambda^{(p)}(G_n)$,
$$x_{\textup{min}}^{p}\geq\frac{1}{n}\Big(1-\frac{p}{2(p-1)\log n}\Big).$$
\end{lem}

\begin{proof}
Suppose, for a contradiction,  that there exists $n_0$ such that for every $n>n_0$, there exist $G_n\in\mathrm{EX}_\lambda^{(p)}(n,F)$ and a principal eigenvector $\mathbf{x}=(x_1,\ldots,x_n)$ of $G_n$ satisfying
\begin{equation}\label{e6}
x_{\textup{min}}^{p}<\frac{1}{n}\Big(1-\frac{p}{2(p-1)\log n}\Big).
\end{equation}

Let  $\lambda_{n}:=\mathrm{ex}_{\lambda}^{(p)}(n,F)$. By Lemma \ref{t3}, we can select $n>n_{0}$ sufficiently large so that
$$\frac{\lambda_{n-1}(n-1)^{2/p-1}}{n-2}-\frac{\lambda_{n}n^{2/p-1}}{n-1}<\frac{1}{n\log n}\cdot\frac{\lambda_{n}n^{2/p-1}}{n-1},$$
which implies
\begin{equation}\label{e7}
\frac{\lambda_{n-1}}{\lambda_{n}}<\frac{n^{2/p-1}(n-2)}{(n-1)^{2/p}}\Big(1+\frac{1}{n\log n}\Big).
\end{equation}
Let $k\in V(G_{n})$ be a vertex with $x_{k}=x_{\textup{min}}$, and let $\mathbf{y}$ be the $(n-1)$-vector obtained from $\mathbf{x}$ by deleting the entry  $x_{k}$. For the subgraph
$G_{n}-k$, the eigenequation \eqref{tzfc} yields
$$ \mathbf{y}^{T}A(G_{n}-k)\mathbf{y}=\lambda^{(p)}(G_{n})-2x_{k}\sum_{i\in N_{G_{n}}(k)}x_{i}=
\lambda_{n}-2\lambda_{n}x_{k}^{p}.$$
Since $G_{n}-k$ is $F$-free, we have
$$\lambda_{n}(1-2x_{k}^{p})=\mathbf{y}^{T}A(G_{n}-k)\mathbf{y}\leq \lambda^{(p)}(G_{n}-k)(\|\mathbf{y} \|_{p})^{2}
\leq \lambda_{n-1}(1-x_{k}^{p})^{2/p}.$$
Thus
\begin{equation}\label{e8}
\frac{\lambda_{n-1}}{\lambda_{n}}\geq \frac{1-2x_{k}^{p}}{(1-x_{k}^{p})^{2/p}}.
\end{equation}

Combining (\ref{e7}) and (\ref{e8}), we obtain
$$\frac{1-2x_{k}^{p}}{(1-x_{k}^{p})^{2/p}}
\leq \frac{n^{2/p-1}(n-2)}{(n-1)^{2/p}}\Big(1+\frac{1}{n\log n}\Big).$$
By Fact $1$ and  (\ref{e6}), we derive
$$\frac{1-\frac{2}{n}(1-\frac{p}{2(p-1)\log n})}{\big({1-\frac{1}{n}(1-\frac{p}{2(p-1)\log n})}\big)^{2/p}}
\leq \frac{1-2x_{k}^{p}}{(1-x_{k}^{p})^{2/p}}\leq\frac{n^{2/p-1}(n-2)}{(n-1)^{2/p}}\Big(1+\frac{1}{n\log n}\Big).$$
Consequently,
$$\frac{\big(n-2+\frac{p}{(p-1)\log n}\big)n^{2/p-1}}
{\big(n-1+\frac{p}{2(p-1)\log n}\big)^{2/p}}
\leq\frac{n^{2/p-1}(n-2)}{(n-1)^{2/p}}\Big(1+\frac{1}{n\log n}\Big).$$
After canceling $n^{2/p-1}$ and rearranging, this becomes
\begin{equation}\label{key-ineq}
1+\frac{p}{(p-1)(n-2)\log n}\leq\Big(1+\frac{p}{2(p-1)(n-1)\log n}\Big)^{2/p}\Big(1+\frac{1}{n\log n}\Big).
\end{equation}

For sufficiently large $n$, we have
\begin{displaymath}
\begin{split}
\Big(1+\frac{p}{2(p-1)(n-1)\log n}\Big)^{2/p}
&=1+\frac{1}{(p-1)(n-1)\log n}+O\Big(\frac{1}{(n\log n)^{2}}\Big)\\
&\leq1+\frac{1}{(p-1)(n-1)\log n}+\frac{1}{(p-1)(n-1)(n-2)\log n}\\
&=1+\frac{1}{(p-1)(n-2)\log n}.
\end{split}
\end{displaymath}
Substituting this estimate into \eqref{key-ineq} yields
$$1+\frac{p}{(p-1)(n-2)\log n}\leq \Big(1+\frac{1}{(p-1)(n-2)\log n}\Big)\Big(1+\frac{1}{n\log n}\Big).$$
After simplification, we get
$$\frac{p}{n-2}\leq \frac{1}{n-2}+\frac{p-1}{n}+\frac{1}{n(n-2)\log n}.$$
which yields $$2(p-1)\leq\frac{1}{\log n}.$$
This is a contradiction for sufficiently large $n$,  completing the proof.
\end{proof}

\begin{lem}\label{p-mindu}
 Let  $p>1$, and let $F$ be a graph with $\pi(F)>0$. Let  $G\in \mathrm{EX}_{\lambda}^{(p)}(n,F)$, and let $\mathbf{x}=(x_{1},\ldots,x_{n})$ be a principal  eigenvector corresponding  to $\lambda^{(p)}(G)$.  Suppose that $0<\varepsilon<1$ and $0<\varepsilon'<\pi(F)\varepsilon$. If  $x_{\textup{min}}^{p}>\frac{1-\varepsilon'}{n}$, then for sufficiently large $n$,
$$\delta(G)>(1-\varepsilon)\pi(F)n.$$
\end{lem}

\begin{proof}
Set  $\lambda:=\lambda^{(p)}(G)$, $\delta:=\delta(G)$, and let $u\in V(G)$ be a vertex attaining minimum degree $\delta(G)$. Considering the eigenequation for $\lambda^{(p)}(G)$  at the vertex  $u$, we have
$$\lambda x_{\textup{min}}^{p-1}\leq \lambda x^{p-1}_{u}=\sum_{v\in N_{G}(u)}x_v.$$
 By H\"older's inequality,
$$(\lambda x_{\textup{min}}^{p-1})^{p}
\leq \delta^{p-1}\sum_{v\in N_{G}(u)}x_v^{p}.$$
Assume, for contradiction, that $\delta(G)\leq(1-\varepsilon)\pi(F)n$.
Since  $x_{\textup{min}}^{p}>\frac{1-\varepsilon'}{n}$ and $0<\varepsilon'<\pi(F)\varepsilon$, we obtain
\begin{displaymath}
\begin{split}
\sum_{v\in N_{G}(u)}x_v^{p}
&= 1-\sum_{w\in V(G)\backslash N_{G}(u)}x_w^{p}\\
&\leq 1-(n-\delta)x_{\textup{min}}^{p}\\
&\leq 1-(1-(1-\varepsilon)\pi(F))(1-\varepsilon')\\
&\leq 1-(1-(1-\varepsilon)\pi(F))(1-\pi(F)\varepsilon)\\
&=  \pi(F)-(1-\varepsilon)\varepsilon\pi(F)^{2}\\
&\leq \pi(F).
\end{split}
\end{displaymath}
Consequently,
$$(\lambda x_{\textup{min}}^{p-1})^{p}
\leq (1-\varepsilon)^{p-1}\pi(F)^{p}n^{p-1}.$$

Choose $0<\varepsilon''<1-(\frac{1-\varepsilon}{1-\varepsilon'})^{\frac{p-1}{p}}$. By Theorem~\ref{p-ess},  for  sufficiently large  $n$,
$$\lambda\geq (1-\varepsilon'')\pi(F)n^{2-2/p}.$$
Using $x_{\mathrm{min}}^p>\frac{1-\varepsilon'}{n}$ again, we obtain
$$(1-\varepsilon'')^{p}\pi(F)^{p}n^{2(p-1)}\cdot \frac{(1-\varepsilon')^{p-1}}{n^{p-1}}
\leq (1-\varepsilon)^{p-1} \pi(F)^{p}n^{p-1},$$
that is
$$1-\varepsilon''\leq \Big(\frac{1-\varepsilon}{1-\varepsilon'}\Big)^{\frac{p-1}{p}},$$
which contradicts the choice of $\varepsilon''$.
\end{proof}

\begin{lem}\label{join}
Let $p>1$,  and let $r,t$ be nonnegative integers with $r\ge 2$.
Write $c=(r-1)/r$ and
\begin{equation}\label{eq:beta}
\beta_{r,p}
   = \frac{2(p-1)}{p}\left[
       \left(\frac{r}{\,r-1\,}\right)^{\!1/(p-1)}
       - \frac{r-1}{r}
     \right].
\end{equation}
Then for sufficiently large $n,$
\begin{equation}\label{eq:main2}
\lambda^{(p)}(K_t\vee T_{n-t,r})
   = cn^{2-2/p} + \beta_{r,p}tn^{1-2/p} + O(n^{-2/p}).
\end{equation}
\end{lem}

\begin{lem}\label{lem:graph-balance}
Let \(p>1\),  
and let
\(r,t, n_1,\ldots,n_r\) be positive integers satisfying  \(r\ge 2\), \(t\ge 1\), $\sum_{i=1}^r n_i=n-t+1$.
Then
\[
\lambda^{(p)}
\!\left(
K_{t-1}\vee K_r(n_1,\ldots,n_r)
\right)
\le
\lambda^{(p)}
\!\left(
K_{t-1}\vee T_{n-t+1,r}
\right).
\]
Moreover, equality holds if and only if $|n_i-n_j|\le 1$ for all $i,j\in[r]$.
\end{lem}

\noindent {\bf Remark.}   The proofs of Lemmas \ref{join} and \ref{lem:graph-balance} are lengthy and we defer  them to the
Appendix.

Firstly, we show that the conclusion of Theorem \ref{main-thm}  holds when the minimum entry of the principal eigenvector corresponding to the $p$-spectral radius of the extremal graph satisfies the quantitative lower
bound in Lemma~\ref{pro-thm}.


\begin{lem}\label{pro-thm}
Let $p>1,$ and let $t, a_1, \ldots, a_r$ be positive integers with $a_1 \ge 2$ and $r \ge 2$.
For sufficiently large $n$, if $G \in \mathrm{EX}_{\lambda}^{(p)}(n, tK^+_r(a_1,\ldots,a_r))$ and
$\mathbf{x} = (x_1, \ldots, x_n)$ is a principal eigenvector corresponding to $\lambda^{(p)}(G)$
satisfying $x_{\min}^p \ge \frac{1}{n}\bigl(1-\frac{p}{2(p-1)\log n}\bigr)$, then
$G \cong K_{t-1} \vee T_{n-t+1,\, r}$.
\end{lem}

\section{Proof of Lemma \ref{pro-thm}}

In this section, we prove Lemma \ref{pro-thm}.  We first recall the concept of degree-stability.
Let $F$ be a graph with $\pi(F) > 0$, and let $\mathscr{H}$ be a family of $F$-free graphs. We say that $F$ is \emph{degree-stable} with respect to $\mathscr{H}$ if there exist constants $n_0$ and $\varepsilon > 0$ such that every $F$-free graph $G$ on $n \ge n_0$ vertices with $\delta(G) \geq (\pi(F) - \varepsilon)n$ belongs to $\mathscr{H}$. Erd\H{o}s and Simonovits \cite{ES1973} proved that color-critical graphs are degree-stable with respect to the family of $r$-partite graphs.

\begin{lem}[\cite{ES1973}]\label{lem-ES1973}
Let $F$ be a color-critical graph with $\chi(F) = r+1 \ge 3$. There exists $n_0$ such that if $G$ is an $F$-free graph on $n \ge n_0$ vertices with $\delta(G) > \frac{3r-4}{3r-1}n$, then $G$ is $r$-partite.
\end{lem}

Now we proceed under the assumption of Lemma~\ref{pro-thm}, i.e., $x_{\min}^p \ge \frac{1}{n}\bigl(1-\frac{p}{2(p-1)\log n}\bigr)$. Let $h = |V(K^+_r(a_1,\ldots,a_r))|$. Fix a  sufficiently small constant $\varepsilon$.  Then $x_{\min}^p \ge \frac{1}{n}(1 - (\frac{\varepsilon(r-1)}{r})^2)$ holds for sufficiently large $n$. By Lemma
~\ref{p-mindu}, we obtain
\begin{equation}\label{p-zxd}
\delta(G)>\frac{(r-1)(1-\varepsilon^{2})}{r}n>\Big(\frac{r-1}{r}-\varepsilon^{2}\Big)n.
\end{equation}
Clearly, $G$ contains a copy of $(t-1)K^{+}_{r}(a_{1},\ldots,a_{r})$.
Otherwise, adding any edge \(e\notin E(G)\) to \(G\) yields a graph that is still \(tK^{+}_{r}(a_{1},\ldots,a_{r})\)-free but has a larger \(p\)-spectral radius, contradicting the assumption that \(G\in \mathrm{EX}_{\lambda}^{(p)}(n,tK^+_r(a_1,\ldots,a_r))\).
Denote such a copy by $H$.  Then $G' = G - V(H)$ is $K^+_r(a_1,\ldots,a_r)$-free. For any vertex $v \in V(G')$,
$$d_{G'}(v)\geq \delta(G)-(t-1)h>\Big(\frac{r-1}{r}-2\varepsilon^{2}\Big)n.$$
Observe that $\frac{r-1}{r}-2\varepsilon^{2}>\frac{3r-4}{3r-1}$. Hence, by Lemma \ref{lem-ES1973}, $G'$ is $r$-partite. Consequently, there exists a vertex partition
$V(G')=V'_{1}\cup\cdots\cup V'_{r}$ such that $\sum_{i=1}^{r}e(G'[V'_{i}])=0$.

 \begin{lem}
There exists a partition $V(G)=V_{1}\cup\cdots\cup V_{r}$ with $V_{i}\supseteq V'_{i}$
for each $i\in[r]$, such that $\sum_{1\leq i<j\leq r}e(V_{i},V_{j})$ is maximum and $||V_{i}|-\frac{n}{r}|<\varepsilon n$ for all $i\in[r]$.
\end{lem}
\begin{proof}
Among  all partitions $V(G)=V_{1}\cup\cdots\cup V_{r}$ with \(V_i \supseteq V_i'\) for each $i\in[r]$, choose  one for which  $\sum_{1\leq i<j\leq r}e(V_{i},V_{j})$ is maximum.
Let $\max_{1\leq i\leq r} \big||V_{i}|-\frac{n}{r}\big|=d$. Without loss of generality assume that $\big||V_{1}|-\frac{n}{r}\big|=d$. Then
\begin{displaymath}
\begin{split}
\frac{n\delta(G)}{2}\leq e(G)&=\sum_{1\leq i<j\leq r}e(V_{i},V_{j})+\sum_{i=1}^{r}e(G[V_{i}])\\
&\leq \sum_{1\leq i<j\leq r}|V_{i}||V_{j}|+(t-1)hn\\
&=|V_{1}|(n-|V_{1}|)+\sum_{2\leq i<j\leq r}|V_{i}||V_{j}|+(t-1)h n.
\end{split}
\end{displaymath}
By the Cauchy-Schwarz inequality, we have $\left(\sum_{i=2}^{r}\left|V_{i}\right|\right)^{2}\leq(r-1)\sum_{i=2}^{r}\left|V_{i}\right|^{2}$, and so
$$2\sum_{2\leq i<j\leq r}|V_{i}||V_{j}|=\Big(\sum_{2\leq i\leq r}|V_{i}|\Big)^{2}-\sum_{2\leq i\leq r}|V_{i}|^{2}
\leq\frac{r-2}{r-1}\Big(\sum_{2\leq i\leq r}|V_{i}|\Big)^{2}=\frac{r-2}{r-1}(n-|V_{1}|)^{2}.$$
Together with \eqref{p-zxd}, we get
\begin{displaymath}
\begin{split}
\Big(\frac{r-1}{2r}-\frac{\varepsilon^{2}}{2}\Big)n^{2}
&\leq|V_{1}|(n-|V_{1}|)+\frac{r-2}{2(r-1)}(n-|V_{1}|)^{2}+(t-1)h n\\
&=\frac{r-1}{2r}n^{2}-\frac{r}{2(r-1)}d^{2}+(t-1)h n,
\end{split}
\end{displaymath}
which yields
$$d\leq\sqrt{\varepsilon^{2}n^{2}-\frac{\varepsilon^{2}}{r}n^{2}+\frac{2(t-1)(r-1)h}{r}n}<\varepsilon n.$$
This completes the proof.
\end{proof}

Let $W=\cup_{i=1}^{r}W_{i}$, where $W_{i}=\{v\in V_{i}: d_{V_{i}}(v)>2th\}$. Since $|V(H)|=(t-1)h$ and  $e(G[V_{i}'])=0$ for each $1\leq i\leq r$, we have
$W\subseteq V(H)$, and hence $|W|\leq th$.

\begin{lem}\label{contains}
For any $S\subseteq V(G)$ with $|S|\leq th$, if there exist a vertex $u\in V_{i}\backslash W$ and an edge $uv\in E(G[V_{i}])$
for some $i\in[r]$, then $G-((W\cup S)\backslash\{u,v\})$ contains a copy of $K^{+}_{r}(a_{1},\ldots,a_{r})$.
\end{lem}
\begin{proof}
Without loss of generality, assume that $u\in V_{1}\backslash W$. By \eqref{p-zxd}, we have
$$\delta(G)>\Big(\frac{r-1}{r}-\varepsilon^{2}\Big)n>\Big(1-\frac{1}{r}-\varepsilon\Big)n.$$
\medskip
\noindent
\textbf{Case 1: $v\in W$.}

Recall that $V(G)=\cup_{i=1}^{r}V_{i}$ is the vertex partition for which $V_{i}\supseteq V'_{i}$
for each $i\in[r]$ and $\sum_{1\leq i<j\leq r}e(V_{i},V_{j})$ is maximum. Hence $d_{V_{1}}(v)\leq\frac{1}{r}d_{G}(v)$. Then for $j\in [r]\backslash\{1\}$, we have

\begin{displaymath}
\begin{split}
d_{V_{j}}(v)
&\geq d_{G}(v)-d_{V_{1}}(v)-(r-2)\Big(\frac{n}{r}+\varepsilon n\Big)\\
&\geq\Big(1-\frac{1}{r}\Big)\Big(1-\frac{1}{r}-\varepsilon\Big)n-(r-2)\Big(\frac{n}{r}+\varepsilon n\Big)\\
&=\frac{n}{r^{2}}-\Big(r-1-\frac{1}{r}\Big)\varepsilon n>\frac{n}{2r^{2}}.
\end{split}
\end{displaymath}
Choose $v_{1,1},\ldots,v_{1,a_{1}-1}\in V_{1}\backslash(W\cup S\cup \{u\})$ and set $u=v_{1,a_{1}}$.
Then for $i\in[a_{1}]$, we have $d_{V_{1}}(v_{1,i})\leq 2t h$, and for $j\in [r]\backslash\{1\}$,
\begin{align}\label{in1}
\begin{split}
d_{V_{j}}(v_{1,i})
&\geq d_{G}(v_{1,i})-(r-2)\Big(\frac{n}{r}+\varepsilon n\Big)-d_{V_{1}}(v_{1,i})\\
&\geq \Big(1-\frac{1}{r}-\varepsilon\Big)n-(r-2)\Big(\frac{n}{r}+\varepsilon n\Big)-2t h\\
&=\frac{n}{r}-(r-1)\varepsilon n-2t h\\
&>\frac{n}{r}-r\varepsilon n.
\end{split}
\end{align}
Then by Lemma \ref{dx},

\begin{displaymath}
\begin{split}
\bigg|\Big(N_{V_{2}}(v)\mathbin{\textstyle\bigcap}\big(\mathbin{\textstyle\bigcap}_{i=1}^{a_{1}}N_{V_{2}}(v_{1,i})\big)\Big)\backslash(W\cup S)\bigg|
&\geq d_{V_{2}}(v)+\sum_{i=1}^{a_{1}}d_{V_{2}}(v_{1,i})-a_{1}|V_{2}|-|W|-|S|\\
&\geq \frac{n}{2r^{2}}+a_{1}\Big(\frac{n}{r}-r\varepsilon n-\frac{n}{r}-\varepsilon n\Big)-2t h\\
&>\frac{n}{2r^{2}}-h(r+2)\varepsilon
n\\
&>a_{2}.
\end{split}
\end{displaymath}

 Choose $a_{2}$ distinct common neighbors
  of $v, v_{1,1}, \ldots,v_{1,a_{1}}$ in $V_{2}\backslash (W\cup S)$, denoted by $v_{2,1},\ldots,v_{2,a_{2}}.$
Now suppose for some integer $k$ with $2\leq k\leq r-1$ that we have already selected vertices $v_{k,1},\ldots,v_{k,a_{k}}$ as  common neighbors of $\{v\}\cup\{v_{i,1},\ldots,v_{i,a_{i}}|1\leq i\leq k-1\}$ in $V_{k}\backslash(W\cup S)$. Similar to (\ref{in1}), we have
$d_{V_{k+1}}(v_{i,j})>\frac{n}{r}-r\varepsilon n$ for $i\in[k]$ and $j\in [a_{i}]$. Note that $\sum_{i=1}^{k}a_{i}\leq h$. Therefore

\begin{displaymath}
\begin{split}
&\bigg|\Big(N_{V_{k+1}}(v)\mathbin{\textstyle\bigcap}\big(\mathbin{\textstyle\bigcap}_{i\in[k], j\in [a_{i}]}N_{V_{k+1}}(v_{i,j})\big)\Big)\backslash(W\cup S)\bigg|\\
&\geq d_{V_{k+1}}(v)+\sum_{i=1}^{k}\sum_{j=1}^{a_{i}}d_{V_{k+1}}(v_{i,j})-\Big(\sum_{i=1}^{k}a_{i}\Big)|V_{k+1}|-|W|-|S|\\
&\geq \frac{n}{2r^{2}}+\sum_{i=1}^{k}a_{i}\Big(\frac{n}{r}-r\varepsilon n-\frac{n}{r}-\varepsilon n\Big)-2t h\\
&>\frac{n}{2r^{2}}-h(r+2)\varepsilon n\\
&>a_{k+1}.
\end{split}
\end{displaymath}

Thus we can select  $v_{k+1,1},\ldots,v_{k+1,a_{k+1}}$ as common neighbors of
$\{v\}\cup\{v_{i,1},\ldots,v_{i,a_{i}}:1\leq i\leq k\}$ in $V_{k+1}\backslash(W\cup S)$. Proceeding inductively,  the subgraph induced by $\{v\}\cup\{v_{i,1},\ldots,v_{i,a_{i}}:1\leq i\leq r\}$ contains a copy of $K^{+}_{r}(a_{1},\ldots,a_{r})$. Consequently, $G-((W\cup S)\backslash\{u,v\})$ contains a $K^{+}_{r}(a_{1},\ldots,a_{r})$.

\medskip
\noindent
\textbf{Case 2: $v\notin W$.}

Note that for $j\in [r]\backslash\{1\}$, $d_{V_{j}}(v)>\frac{n}{r}-r\varepsilon n$.
Arguing as in Case~1, we obtain  that $G-((W\cup S)\backslash\{u,v\})$ contains a $K^{+}_{r}(a_{1},\ldots,a_{r})$.
\end{proof}

Let $x_{\textup{max}}=\max\{x_{1},\ldots,x_{n}\}$, and let $v^{*}$ be a vertex in $V(G)\backslash W$ such that $x_{v^{*}}=\max_{v\in V(G)\backslash W}x_{v}$.

\begin{lem}\label{up-bound}
$x_{v^{*}}\leq 4^{\frac{1}{p-1}}n^{-\frac{1}{p}}$.
\end{lem}

\begin{proof}
Since $K_{t-1}\vee T_{n-t+1,r}$ is $tK^+_r(a_1,\ldots,a_r)$-free, Lemma
 \ref{join} yields $$\lambda^{(p)}(G)\geq\Big(1-\frac{1}{r}\Big)n^{\frac{2(p-1)}{p}}+O(n^{-2/p})=\Big(1-\frac{1}{r}+o(1)\Big)n^{\frac{2(p-1)}{p}}.$$
By $x_{\min}^p \ge \frac{1}{n}(1-\frac{p}{2(p-1)\log n})$, we have
$$\lambda^{(p)}(G)x_{v^{*}}^{p-1}\geq \Big(1-\frac{1}{r}+o(1)\Big)n^{\frac{2(p-1)}{p}}\cdot \Big(\frac{1}{n}\Big(1-\frac{p}{2(p-1)\log n}\Big)\Big)^{\frac{p-1}{p}}\geq \frac{r-1}{2r}n^{\frac{p-1}{p}}.$$

Assume that $v^{*}\in V_{i_{0}}$, where $1\leq i_{0}\leq r$. Then $|N_{V_{i_{0}}\backslash W}(v^{*})|\leq|N_{V_{i_{0}}}(v^{*})|\leq 2th$ as $v^{*}\notin W$. Since $|W|\leq th$, we have
\begin{align*}
\lambda^{(p)}(G)x_{v^{*}}^{p-1}
&= \biggl(\sum_{v\in N_{W}(v^{*})}x_{v}\biggr)
   + \biggl(\sum_{v\in N_{V_{i_{0}}\backslash W}(v^{*})}x_{v}\biggr)
   + \biggl(\sum_{v\in \cup_{i\neq i_{0}}N_{V_{i}\backslash W}(v^{*})}x_{v}\biggr)\\
&\leq th x_{\textup{max}} + 2th x_{v^{*}} + \sum_{v\in\cup_{i\neq i_{0}}V_{i}\backslash W}x_{v}\\
&\leq \frac{1}{3}\lambda^{(p)}(G)x_{v^{*}}^{p-1} + \sum_{v\in\cup_{i\neq i_{0}}V_{i}\backslash W}x_{v}.
\end{align*}
 By $r\geq2$ and H\"older's inequality, we obtain
\begin{align*}
\frac{1}{4}n^{\frac{2(p-1)}{p}}x_{v^{*}}^{p-1}
&\leq \frac{2}{3}\lambda^{(p)}(G)x_{v^{*}}^{p-1}\\
&\leq  \sum_{v\in\cup_{i\neq i_{0}}V_{i}\backslash W}x_{v}\\
&\leq \bigg(\sum_{v\in\cup_{i\neq i_{0}}V_{i}\backslash W}x_{v}^{p}\bigg)^{\frac{1}{p}}
(|V(G)|)^{\frac{p-1}{p}}\\
&\leq n^{\frac{p-1}{p}},
\end{align*}
which implies $x_{v^{*}}\leq 4^{\frac{1}{p-1}}n^{-\frac{1}{p}}$. This completes the proof.
\end{proof}

\begin{lem}\label{W}
$|W|=t-1$.
\end{lem}
\begin{proof}
When $t=1,$ we have $|W|=0$ by $W\subseteq V(H)=\emptyset,$  which is  the desired equality. In the sequel, assume that $t\geq 2.$
We  first claim  that $|W|\leq t-1$. Suppose otherwise that $|W|\geq t$. Let $S_{1}=\emptyset$
and $\{w_{1},w_{2},\ldots,w_{t}\}\subseteq W$, where $w_{j}\in V_{i_{j}}$ for each $j\in [t]$.
Since $|W|\leq th$ and $d_{V_{i_1}}(w_{1})>2th$, there exists $v_{1}\in V_{i_1}\backslash (W\cup S_{1})$ such that $w_{1}v_{1}\in E(G[V_{i_1}])$. By Lemma \ref{contains},
$G-((W\cup S_{1})\backslash\{w_{1},v_{1}\})$ contains a $K^{+}_{r}(a_{1},\ldots,a_{r})$, say $H_{1}$. Set $S_{2}=S_{1}\cup V(H_{1})$.
Note that  $|W\cup S_{2}|\leq 2t h$, so there exists $v_{2}\in V_{i_2}\backslash (W\cup S_{2})$ such that
$w_{2}v_{2}\in E(G[V_{i_2}])$. Since $|S_{2}|\leq h$, applying  Lemma \ref{contains},
$G-((W\cup S_{2})\backslash\{w_{2},v_{2}\})$ contains a $K^{+}_{r}(a_{1},\ldots,a_{r})$, say $H_{2}$.
For each $3\leq k\leq t$, let $S_{k}=S_{k-1}\cup V(H_{k-1}).$
Then $|S_{k}|=(k-1)h\leq th$ and hence $|W\cup S_{k}|\leq 2t h.$
So there exists $v_{k}\in V_{i_k}\backslash (W\cup S_{k})$ such that
$w_{k}v_{k}\in E(G[V_{i_k}])$.
Applying  Lemma \ref{contains} again,
$G-((W\cup S_{k})\backslash\{w_{k},v_{k}\})$ contains a $K^{+}_{r}(a_{1},\ldots,a_{r})$, say $H_{k}$.
Consequently, $H_{1},\ldots,H_{t}$ are $t$ vertex-disjoint copies of $K^{+}_{r}(a_{1},\ldots,a_{r})$ in $G$, a contradiction.

Thus, $|W|\leq t-1$. 
Now suppose  that  $|W|\leq t-2$. Recall that $H$ is a copy of $(t-1)K^{+}_{r}(a_{1},\ldots,a_{r})$ in $G$, $W\subseteq V(H)$, and $G-V(H)$ is $r$-partite. Define $G''=\cup_{i=1}^{r}G[V_{i}\backslash W]$. Then $e(G'')\leq 2t^{2}h^{2}$. Choose $u_{0}\in V_{1}\backslash W$ and construct
$$G^{*}=G- E(G'') +\{u_{0}v_{0}:v_{0}\in V_{1}\backslash (W\cup\{u_{0}\})\}.$$
Then $G^{*}$ is $tK^{+}_{r}(a_{1},\ldots,a_{r})$-free.
From $x_{\min}^p \ge \frac{1}{n}(1-\frac{p}{2(p-1)\log n})$ and Lemma \ref{up-bound}, we have
\begin{displaymath}
\begin{split}
\lambda^{(p)}(G^{*})-\lambda^{(p)}(G)
&\geq 2\sum_{v_{0}\in V_{1}\backslash (W\cup\{u_{0}\})}x_{u_{0}}x_{v_{0}}-2\sum_{uv\in E(G'')}x_{u}x_{v}\\
&\geq 2n^{-\frac{2}{p}}\Big(1-\frac{p}{2(p-1)\log n}\Big)^{\frac{2}{p}}\Big(\frac{n}{r}-2\varepsilon n\Big)
-4t^{2}h^{2}\cdot 4^{\frac{2}{p-1}}n^{-\frac{2}{p}}\\
&\geq \Big(\frac{1}{r}-2\varepsilon\Big)n^{1-\frac{2}{p}}-4^{\frac{2}{p-1}+1}t^{2}h^{2}n^{-\frac{2}{p}}>0,
\end{split}
\end{displaymath}
which contradicts the maximality of $\lambda^{(p)}(G)$. Therefore $|W|=t-1$, completing the proof.
\end{proof}

\begin{proof}[{\bf Proof of Lemma \ref{pro-thm}}]
We first claim that $e(G[V_{i}\backslash W])= 0$ for each $1\leq i\leq r.$ Otherwise,
if $e(G[V_{i}\backslash W])\neq 0$ for some $1\leq i\leq r,$ then, by applying a construction similar to that in Lemma~\ref{W},  we can obtain $t$ vertex-disjoint copies of $K^{+}_{r}(a_{1},\ldots,a_{r})$ in $G$ using Lemmas \ref{contains} and \ref{W},  a contradiction.

Hence, $G-W$ is an $r$-partite graph. Let $n_{i}=|V_{i}\backslash W|$ for each $1\leq i\leq r.$ Then
$\sum_{i=1}^{r}n_{i}=n-t+1$ and $G$ is a copy of a spanning subgraph of $K_{t-1}\vee K_{r}(n_{1},\ldots, n_{r})$.
Noting that $x_{\min}>0,$ we have $G \cong K_{t-1}\vee K_{r}(n_{1},\ldots, n_{r}).$
Combining with
 Lemma \ref{lem:graph-balance},  $G \cong K_{t-1} \vee T_{n-t+1,\, r}$.
\end{proof}

\section{Proof of Theorem \ref{main-thm}}

\begin{proof}[{\bf Proof of Theorem \ref{main-thm}}]
By Lemmas \ref{t4} and \ref{pro-thm}, there exists a sufficiently large $n_{0}$ such that
$$\mathrm{EX}_{\lambda}^{(p)}(n_{0},tK^{+}_{r}(a_{1},\ldots,a_{r}))=\{K_{t-1}\vee T_{n_{0}-t+1,r}\}.$$
We proceed by induction on $n$. The case $n=n_{0}$ is the base above.
Assume the statement holds for $n-1\ge n_{0}$; we prove it for $n$.

Let $G\in \mathrm{EX}_{\lambda}^{(p)}(n,tK^{+}_{r}(a_{1},\ldots,a_{r}))$ and $\mathbf{x}=(x_{1},\ldots,x_{n})$ be a principal eigenvector corresponding to $\lambda^{(p)}(G)$. By Lemma \ref{pro-thm}, it suffices to show that
$x_{\min}^p \ge \frac{1}{n}(1-\frac{p}{2(p-1)\log n})$.
Set $F=tK^{+}_{r}(a_{1},\ldots,a_{r})$.
Since $(1-x)^{-2/p}=1+\frac{2}{p}x+O(x^{2})$ as $x\to 0$, and since
$x_{\min}^{p}\le 1/n$, inequality \eqref{e8} yields
\begin{equation}\label{in2}
\frac{\mathrm{ex}_{\lambda}^{(p)}(n-1,F)}{\mathrm{ex}_{\lambda}^{(p)}(n,F)}\geq \frac{1-2x_{\min}^{p}}{(1-x_{\min}^{p})^{2/p}}
=1-2\Big(1-\frac{1}{p}\Big)x_{\min}^{p}+O(n^{-2}).
\end{equation}

On the other hand, since $K_{t-1}\vee T_{n-t+1,r}$ is $F$-free, it follows
from Lemma \ref{join} and the induction hypothesis that
\begin{align}\label{in3}
\begin{split}
\frac{\mathrm{ex}_{\lambda}^{(p)}(n-1,F)}{\mathrm{ex}_{\lambda}^{(p)}(n,F)}
&=\frac{\lambda^{(p)}(K_{t-1}\vee T_{n-t,r})}{\mathrm{ex}_{\lambda}^{(p)}(n,F)}\\
&\leq \frac{\lambda^{(p)}(K_{t-1}\vee T_{n-t,r})}{\lambda^{(p)}(K_{t-1}\vee T_{n-t+1,r})}\\
&=\frac{c(n-1)^{2-2/p} + \beta_{r,p}(t-1)(n-1)^{1-2/p} + O(n^{-2/p})}{cn^{2-2/p} + \beta_{r,p}(t-1)n^{1-2/p} + O(n^{-2/p})}\\
&=1-\frac{2c(1-1/p)n^{1-2/p}}{cn^{2-2/p} + \beta_{r,p}(t-1)n^{1-2/p} + O(n^{-2/p})}+O(n^{-2})\\
&\leq1-\frac{2c(1-1/p)}{cn + \beta_{r,p}(t-1)+1}+O(n^{-2})=1-\frac{2(1-1/p)}{n}+O(n^{-2}),
\end{split}
\end{align}
where $c$ and $\beta_{r,p}$ are defined in Lemma  \ref{join}.

Combining (\ref{in2}) and  (\ref{in3}) yields
$$x_{\textup{min}}^{p}\geq\frac{1}{n}(1-O(n^{-1}))\geq\frac{1}{n}\Big(1-\frac{p}{2(p-1)\log n}\Big).$$
This completes the induction and hence the proof of Theorem \ref{main-thm}.
\end{proof}

\subsection*{Acknowledgement}
 The authors developed the core ideas and proof strategy, and acknowledge the use of AI assistance only in discussions concerning the proofs of Lemmas \ref{join} and \ref{lem:graph-balance}.

\section*{Appendix}

\begin{proof}[{\bf Proof of Lemma \ref{join}}]
When $t=0,$ on the one hand, let $\mathbf{x}$ be  the $n$-vector $\mathbf{x}=(n^{-1/p},\ldots,n^{-1/p}).$ By~(\ref{def}), we have $\lambda^{(p)}(T_{n,r})\geq \frac{r-1}{r}n^{2-2/p}-O(n^{-2/p}).$
On the other hand,  we have $\lambda^{(p)}(T_{n,r})\leq  \frac{r-1}{r}n^{2-2/p}$ by~\cite[Theorem~3]{KNY2014}. Thus, $\lambda^{(p)}(T_{n,r})= \frac{r-1}{r}n^{2-2/p}+O(n^{-2/p}).$
In the sequel, assume $t\geq 1.$
Set
$$T_{n-t,r}=K_{a_1,\dots,a_r},\qquad
\sum_{i=1}^{r}a_{i}=n-t,\qquad
|a_{i}-a_{j}|\leq 1.$$
Thus we may write
\begin{equation}\label{eq:ai}
a_i=\frac{n}{r}+\delta_i,\qquad \delta_i=O(1),\qquad
\sum_{i=1}^r\delta_i=-t,
\end{equation}

Let $G_{n}=K_t\vee T_{n-t,r}$, and  let $\mathbf{x}=(x_{1},\ldots,x_{n})$ be a principal  eigenvector corresponding  to $\lambda^{(p)}(G_{n})$.
Then $\sum_{v\in V(G_{n})}x_{v}^{p}=1$. Set $\lambda:=\lambda^{(p)}(G_n)$. By  the eigenequations (\ref{tzfc}), we can write $x_i$ for the common value in part $i$ of $K_{a_1,\dots,a_r},$ $x_0$ for the common value on $K_t$  and thus,
\begin{gather}
\lambda x_0^{p-1}
   = (t-1)x_0+\sum_{i=1}^r a_i x_i, \label{eq:sys0}\\
\lambda x_i^{p-1}
   = t x_0+\sum_{j\ne i}a_j x_j,
      \qquad i=1,\dots,r, \label{eq:sysi}
\end{gather}
together with the normalization
\begin{equation}\label{eq:norm}
t x_0^p+\sum_{i=1}^r a_i x_i^p=1.
\end{equation}
From (\ref{p-bds}), we have
$$2e(G_n)n^{-2/p}\leq  \lambda\leq (2e(G_n))^{1-1/p},$$
which implies
$cn^{2-2/p}+O(n^{1-2/p})\leq \lambda\leq n^{2-2/p}$.

Now introduce the scaled quantities
\begin{equation}\label{eq:scale}
L_n=\frac{\lambda}{n^{2-2/p}},\qquad
x_0=n^{-1/p}y_0,\qquad
x_i=n^{-1/p}y_i.
\end{equation}
Dividing \eqref{eq:sys0}--\eqref{eq:sysi} by $n^{1-1/p}$ and
\eqref{eq:norm} by $1$ gives
\begin{gather}
L_n y_0^{p-1}
   = \frac{t-1}{n}y_0+\sum_{i=1}^r\frac{a_i}{n}\,y_i,
      \label{eq:sc0}\\
L_n y_i^{p-1}
   = \frac{t}{n}y_0+\sum_{j\ne i}\frac{a_j}{n}\,y_j,
      \qquad i=1,\dots,r, \label{eq:sci}\\
\frac{t}{n}y_0^p+\sum_{i=1}^r\frac{a_i}{n}\,y_i^p
   = 1. \label{eq:scnorm}
\end{gather}


By \eqref{eq:scnorm} and $a_i/n\to 1/r$, the $y_i$ are bounded.
Since  $\frac{t-1}{n}y_0=\frac{t-1}{n^{1-1/p}}x_0$, $\frac{t-1}{n}y_0$ is bounded.
Then \eqref{eq:sc0} and the positive lower bound on $L_n$ bound $y_0$.
By the Bolzano--Weierstrass theorem and a standard diagonal argument, we may extract a subsequence $n_{k}$
 such that all components converge simultaneously:
\[
L_{n_{k}}\to L,\qquad y_0\to y_0^*,\qquad y_i\to y_i^*
 \qquad (k\to\infty).
\]
Letting $k\to\infty$ (so $a_i/n_{k}\to 1/r$) in
\eqref{eq:sc0}--\eqref{eq:scnorm} yields
\begin{align}
L\,(y_0^*)^{p-1} &= \frac{1}{r}\sum_{j=1}^r y_j^*, \label{eq:lim0}\\[2pt]
L\,(y_i^*)^{p-1} &= \frac{1}{r}\sum_{j\ne i}y_j^*,
   \qquad i=1,\dots,r, \label{eq:limi}\\[2pt]
\frac{1}{r}\sum_{i=1}^r (y_i^*)^p &= 1. \label{eq:limnorm}
\end{align}
Set $M=\dfrac{1}{r}\sum_{j=1}^r y_j^*$. Then \eqref{eq:limi} can be
written as
\[
L\,(y_i^*)^{p-1}+\frac{1}{r}y_i^*=M.
\]
Since the function $f(s)=Ls^{p-1}+s/r$ is strictly increasing on $s>0$, we have
$y_1^*=\cdots=y_r^*$. From \eqref{eq:limnorm} it follows that
$y_1^*=\cdots=y_r^*=1$; substituting into \eqref{eq:limi} gives
$L=(r-1)/r=c$. Finally \eqref{eq:lim0} gives $c\,(y_0^*)^{p-1}=1$, which implies
\begin{equation}\label{eq:zdef}
y_0^*=z:=c^{-1/(p-1)}
   =\left(\frac{r}{\,r-1\,}\right)^{\!1/(p-1)}.
\end{equation}
Since every convergent subsequence has the same limit, the full
sequences converge:
\begin{equation}\label{eq:limits}
L_n\to c,\qquad y_i\to 1,\qquad y_0\to z \qquad (n\to\infty).
\end{equation}


We now show that
\[
y_i=1+\frac{s_i}{n}\quad \text{with } s_i=O(1),
\qquad
y_0=z+O(n^{-1})
\]
as $n \to\infty$.

Subtracting \eqref{eq:sci} for distinct indices $i$ and $j$ gives
\[
L_n\!\left(y_i^{p-1}-y_j^{p-1}\right)
   = \frac{a_j}{n}y_j-\frac{a_i}{n}y_i.
\]
The mean-value
theorem provides $\xi_{ij}$ between $y_i$ and $y_j$ such that
$y_i^{p-1}-y_j^{p-1}=(p-1)\xi_{ij}^{p-2}(y_i-y_j)$. Combining with  $a_i/n=\tfrac{1}{r}+O(n^{-1})$ by \eqref{eq:ai}, we have
\[
\left[L_n(p-1)\xi_{ij}^{p-2}+\frac{1}{r}\right](y_i-y_j)
   = O(n^{-1}).
\]
The bracket tends to $c(p-1)+\tfrac{1}{r}>0$, so $y_i-y_j=O(n^{-1})$. Consequently, there exists a number $y$ such that
$y_i=y+O(n^{-1})$ for every $i\in[r]$. By the mean value theorem again, there exists $\xi_i$ between $y_i$ and $y$ such that
$$y_i^p-y^p=p\xi_i^{p-1}(y_i-y)=O(n^{-1}),$$
that is, $y_i^p=y^p+O(n^{-1})$. Substituting into \eqref{eq:scnorm} and using $y_0=O(1)$ gives
$y^p=1+O(n^{-1})$. Thus  $y=1+O(n^{-1})$, and hence
\begin{equation}\label{eq:yi_O1n}
y_i=1+O(n^{-1})\qquad (i=1,\dots,r).
\end{equation}
From \eqref{eq:sci}, noting that
$\sum_{j\ne i}\tfrac{a_j}{n}y_j=c+O(n^{-1})$ and
$y_i^{p-1}=1+O(n^{-1})$, we obtain $L_n=c+O(n^{-1})$.
Then \eqref{eq:sc0} gives  $L_n y_0^{p-1}=1+O(n^{-1})$, so
$$y_0=\big(c^{-1}+O(n^{-1})\big)^{1/(p-1)}=z+O(n^{-1}).$$
Consequently we may write
\begin{equation}\label{eq:expansions}
y_i=1+\frac{s_i}{n},\qquad s_i=O(1),
\qquad
y_0=z+O(n^{-1}).
\end{equation}

By \eqref{eq:ai} and  $(1+\frac{s_i}{n})^p=1+\frac{p s_i}{n}+O(n^{-2})$, we have
\begin{align*}
\frac{a_i}{n}\,y_i^p
   &=\left(\frac{1}{r}+\frac{\delta_i}{n}\right)
     \left(1+\frac{p s_i}{n}+O(n^{-2})\right)\\[2pt]
   &= \frac{1}{r}
      +\frac{1}{n}\left(\delta_i+\frac{p}{r}s_i\right)
      +O(n^{-2}).
\end{align*}
Summing over $i$ and using $\sum_{i=1}^{r}\delta_i=-t$ gives
\[
\sum_{i=1}^r\frac{a_i}{n}\,y_i^p
   = 1+\frac{1}{n}\left(-t+\frac{p}{r}\sum_{i=1}^r s_i\right)
     +O(n^{-2}).
\]
By $y_0=z+O(n^{-1})$ in \eqref{eq:expansions}, we have
$\tfrac{t}{n}y_0^p=\tfrac{t z^p}{n}+O(n^{-2})$.
Returning to \eqref{eq:scnorm}, the constant terms cancel and we obtain
\[
-t+\frac{p}{r}\sum_{i=1}^r s_i+t z^p=O(n^{-1}),
\]
that is,
\begin{equation}\label{eq:sbar}
\bar{s}:=\frac{1}{r}\sum_{i=1}^r s_i
   = -\frac{t}{p}\left(z^p-1\right)+O(n^{-1}).
\end{equation}

We now compute $\lambda$ by inserting the expansions
for $x_0$ and $x_i$ into $\lambda=\mathbf{x}^{T}A(G_n)\mathbf{x}$. The edges of $G_n$ are
split into three classes: (i) edges between distinct Tur\'an parts, (ii) join-edges
between $K_t$ and $T_{n-t,r}$, and (iii) edges inside $K_t$. We consider the contribution of each class
in turn.

\medskip
\noindent
\textbf{Case 1: Tur\'an part.}

For convenience, set $B_i=a_i x_i$. By \eqref{eq:ai} and \eqref{eq:expansions},
\[
B_i
   =\frac{1}{r}n^{1-1/p}
    +\left(\delta_i+\frac{s_i}{r}\right)n^{-1/p}
    +O(n^{-1-1/p}).
\]
Summing,
\begin{equation}\label{eq:Bi}
\sum_{i=1}^r B_i
   = n^{1-1/p}
     +(-t+\bar{s})n^{-1/p}
     +O(n^{-1-1/p}),
\end{equation}
and
\[
\sum_{i=1}^r B_i^2
   = \frac{1}{r}\,n^{2-2/p}
     +\frac{2}{r}(-t+\bar{s})n^{1-2/p}
     +O(n^{-2/p}).
\]
The contribution of the Tur\'an graph is therefore
\begin{align}
Q_T:
   &= \Bigl(\sum_{i=1}^{r} B_i\Bigr)^2-\sum_{i=1}^{r} B_i^2 \notag\\[2pt]
   &= cn^{2-2/p}
      +2c(-t+\bar{s})n^{1-2/p}
      +O(n^{-2/p}), \label{eq:QT}
\end{align}
where $c=1-1/r$.

\medskip
\noindent
\textbf{Case 2: Join part.}

Each vertex of $K_t$ is adjacent to every vertex of the Tur\'an part,
so by \eqref{eq:expansions} and \eqref{eq:Bi},
\begin{equation}\label{eq:QJ}
Q_J:=2tx_0\sum_{i=1}^r a_i x_i
   =2t zn^{1-2/p}+O(n^{-2/p}).
\end{equation}

\medskip
\noindent
\textbf{Case 3: Clique part.}

Since $K_t$ has $\binom{t}{2}$ edges and $t$ is fixed,
\begin{equation}\label{eq:QK}
Q_K:=2\tbinom{t}{2}x_0^2
   = O(n^{-2/p}).
\end{equation}

Summing \eqref{eq:QT}, \eqref{eq:QJ}, and \eqref{eq:QK}, we obtain
\[
\lambda
   = cn^{2-2/p}
     +2\left[c(-t+\bar{s})+t z\right]n^{1-2/p}
     +O(n^{-2/p}).
\]
Substituting \eqref{eq:sbar}, and noting that
$c(-t+\bar{s})
   =-ct-\tfrac{ct}{p}(z^p-1)+O(n^{-1})$, we get
\begin{align*}
\lambda
   &= cn^{2-2/p}
      +2t\left[z-c-\frac{c}{p}(z^p-1)\right]n^{1-2/p}
      +O(n^{-2/p}).
\end{align*}
Since $c z^{p-1}=1$ by \eqref{eq:zdef}, multiplication by $z$ gives
$c z^p=z$. Hence
\[
z-c-\frac{c}{p}(z^p-1)
   = z-c-\frac{1}{p}(z-c)
   = \frac{p-1}{p}(z-c).
\]
We conclude
\begin{equation}\label{eq:final}
\lambda^{(p)}(G_n)
   = cn^{2-2/p}
     +\frac{2(p-1)}{p}(z-c)tn^{1-2/p}
     +O(n^{-2/p}).
\end{equation}
Inserting $z=\left(\dfrac{r}{r-1}\right)^{1/(p-1)}$ and
$c=\dfrac{r-1}{r}$ establishes \eqref{eq:main2} and identifies
$\beta_{r,p}$ as in \eqref{eq:beta}.
\end{proof}

\begin{proof} [{\bf Proof of Lemma \ref{lem:graph-balance}}]

We first observe that $K_{t-1}\vee K_r(n_1,\ldots,n_r)$
is a complete multipartite graph with \(t-1\) parts of size $1$ and  $r$  parts of size $n_{1},\ldots,n_{r}$.
Among all such graphs with \(\sum_{i=1}^r n_i=n-t+1\), choose  \(G\) with maximum
\(p\)-spectral radius. We show that the \(r\)  parts of sizes $n_{1},\ldots,n_{r}$  are  as balanced as
 possible.

Suppose, to the contrary, that two of these parts, say \(A\) and \(B\),
have sizes $|A|=a$ and  $|B|=b$, where $b-a\ge 2$. Since $G$ is connected, there clearly exists  a Perron--Frobenius eigenvector
 $\mathbf{x}=(x_{1},\ldots,x_{n})$  corresponding to \(\lambda^{(p)}(G)\).
As vertices in the same part of a complete multipartite graph
have identical neighborhoods, the eigenequations imply that the coordinates
of $\mathbf{x}$ are constant on each part. Write $x_v=\alpha$ for $v\in A$, and $x_v=\beta$ for $v\in B$.

Set
\[
c=a\alpha^p+b\beta^p
\qquad\text{and}\qquad
R=\sum_{v\notin A\cup B}x_v.
\]
Let $C$ be the sum of $2x_{u}x_{v}$ over all edges $uv\in E(G)$
with both endpoints outside \(A\cup B\). Then
\begin{equation}\label{eq:PG}
\mathbf{x}^{T}A(G)\mathbf{x}
=C+2R(a\alpha+b\beta)+2ab\alpha\beta.
\end{equation}

We now balance the two parts \(A\) and \(B\).

\medskip
\noindent
\textbf{Case 1: \(a+b=2m\).}

Replace \(A\) and \(B\) by two parts of sizes \(m\) and \(m\), and denote
the resulting graph by \(G'\). On these two new parts assign the common
value
\[
d=\left(\frac{c}{2m}\right)^{1/p}=\left(\frac{a\alpha^p+b\beta^p}{2m}\right)^{1/p},
\]
while keeping all other coordinates unchanged. Let the resulting vector
be \(\mathbf{y}\). Then \(\|\mathbf{y}\|_p=1\).

By H\"older's inequality,
\begin{equation}\label{eq:jc-bound}
a\alpha+b\beta
\le
(a+b)^{1-1/p}(a\alpha^p+b\beta^p)^{1/p}
=2md.
\end{equation}
Also, by the arithmetic--geometric mean inequality,
\[
a\alpha^p+b\beta^p
\ge
2\sqrt{ab}\,(\alpha\beta)^{p/2}.
\]
Hence
\[
m^2d^2
\ge
m^{2-2/p}(ab)^{1/p}\alpha\beta
=
ab\alpha\beta
\left(\frac{m^2}{ab}\right)^{1-1/p}.
\]
Since \(a+b=2m\) and \(a\ne b\), we have \(ab<m^2\).
As \(p>1\), this yields $m^2d^2>ab\alpha\beta$.
Together with \eqref{eq:PG} and \eqref{eq:jc-bound}, we get
\[
\mathbf{y}^{T}A(G')\mathbf{y}>\mathbf{x}^{T}A(G)\mathbf{x}=\lambda^{(p)}(G),
\]
and therefore $\lambda^{(p)}(G')>\lambda^{(p)}(G)$, contradicting  the maximality of \(G\).

\medskip
\noindent
\textbf{Case 2: \(a+b=2m+1\).}

Since \(b-a\ge 2\) and \(a+b\) is odd, in fact $b-a\ge 3$.
Thus
\begin{equation}\label{eq:ab-bound}
ab\le (m-1)(m+2)<m(m+1).
\end{equation}

Replace \(A\) and \(B\) by two parts of sizes \(m\) and \(m+1\), and
denote the resulting graph again by \(G'\). First assign the common value
\[
d=\left(\frac{c}{2m+1}\right)^{1/p}=\left(\frac{a\alpha^p+b\beta^p}{2m+1}\right)^{1/p}
\]
to all vertices in these two new parts, keeping all other coordinates
unchanged. Let the resulting vector be \(\mathbf{y}\). Again \(\|\mathbf{y}\|_p=1\), and
H\"older's inequality gives
\begin{equation}\label{eq:linear}
a\alpha+b\beta\le
(a+b)^{1-1/p}(a\alpha^p+b\beta^p)^{1/p}=(2m+1)d.
\end{equation}

We distinguish two ranges of \(p\).

\smallskip
\noindent
\emph{Subcase 2.1: \(p>9/8\).}

By the arithmetic--geometric mean inequality,
\[
d^2
\ge
\left(\frac{4ab}{(2m+1)^2}\right)^{1/p}\alpha\beta.
\]
Therefore
\[
\frac{m(m+1)d^2}{ab\alpha\beta}
\ge
\frac{m(m+1)}{ab}
\left(\frac{4ab}{(2m+1)^2}\right)^{1/p}.
\]
The right-hand side is decreasing as a function of \(ab\), since \(1/p-1<0\). We claim that $\frac{m(m+1)d^2}{ab\alpha\beta}>1$.
Using \eqref{eq:ab-bound}, it suffices to prove
\[
\frac{m(m+1)}{(m-1)(m+2)}
>
\left(
\frac{(2m+1)^2}{4(m-1)(m+2)}
\right)^{1/p}.
\]
Now
\[
\frac{m(m+1)}{(m-1)(m+2)}
=
1+\frac{2}{(m-1)(m+2)}
\]
and
\[
\frac{(2m+1)^2}{4(m-1)(m+2)}
=
1+\frac{9}{4(m-1)(m+2)}.
\]
Since \(p>9/8\), we have \(1/p<8/9\). By Bernoulli's inequality,
\[
\left(
1+\frac{9}{4(m-1)(m+2)}
\right)^{1/p}
<
1+\frac{2}{(m-1)(m+2)}.
\]
Hence
\begin{equation}\label{eq:cross-strict}
m(m+1)d^2>ab\alpha\beta.
\end{equation}
Combining \eqref{eq:linear} and \eqref{eq:cross-strict} with
\eqref{eq:PG}, we again obtain
\[
\mathbf{y}^{T}A(G')\mathbf{y}>\mathbf{x}^{T}A(G)\mathbf{x}=\lambda^{(p)}(G),
\]
a contradiction.

\smallskip
\noindent
\emph{Subcase 2.2: \(1<p\le 9/8\).}

If $m(m+1)d^2>ab\alpha\beta$, then the same argument as above gives a contradiction. Thus suppose
instead that
\begin{equation}\label{eq:bad-cross}
m(m+1)d^2\le ab\alpha\beta.
\end{equation}

We now use a different test vector. Define
\[
u=\frac{a\alpha}{m},
\qquad
v=\frac{b\beta}{m+1}.
\]
Assign the value \(u\) to the part of size \(m\), the value \(v\) to the
part of size \(m+1\), and keep all remaining coordinates unchanged.
Let the resulting vector be \(\mathbf{z}\). Then
\[
mu=a\alpha,
\qquad
(m+1)v=b\beta,
\]
and $m(m+1)uv=ab\alpha\beta$. Consequently,
\[
\mathbf{z}^{T}A(G')\mathbf{z}=\mathbf{x}^{T}A(G)\mathbf{x}=\lambda^{(p)}(G).
\]
By the maximality of \(G\), we must have $\|\mathbf{z}\|_p\ge 1$.
Set
\[
\mu=a\alpha^p,
\qquad
\nu=b\beta^p.
\]
The inequality $\|\mathbf{z}\|_p\ge 1$ yields
\[
\mu\left(\frac{a}{m}\right)^{p-1}
+
\nu\left(\frac{b}{m+1}\right)^{p-1}
\ge \mu+\nu.
\]
Since \(a+b=2m+1\), write $a=m-s$ and $b=m+1+s$ for some integer \(s\ge 1\). Hence
\begin{equation}\label{eq:AB}
\nu\left[
\left(1+\frac{s}{m+1}\right)^{p-1}-1
\right]
\ge
\mu\left[
1-\left(1-\frac{s}{m}\right)^{p-1}
\right].
\end{equation}

Because \(0<p-1\le 1/8<1\), the function \(x^{p-1}\) is strictly
concave on \((0,\infty)\). Therefore,
\[
1-\left(1-\frac{s}{m}\right)^{p-1}
>
(p-1)\frac{s}{m},
\]
whereas
\[
\left(1+\frac{s}{m+1}\right)^{p-1}-1
<
(p-1)\frac{s}{m+1}.
\]
Using these inequalities in \eqref{eq:AB}, we obtain
\begin{equation}\label{eq:BoverA}
\frac{\nu}{\mu}>\frac{m+1}{m}.
\end{equation}

Since the function $\frac{x}{(1+x)^2}$
is strictly decreasing for \(x>1\), \eqref{eq:BoverA} implies
\[
\frac{\mu\nu}{(\mu+\nu)^2}
<
\frac{m(m+1)}{(2m+1)^2}.
\]
Recalling that
$\mu\nu=ab(\alpha\beta)^p$ and $\mu+\nu=c$, we obtain
\[
ab(\alpha\beta)^p
<
\frac{m(m+1)}{(2m+1)^2}c^2.
\]
Since $d^2=(\frac{c}{2m+1})^{2/p}$,
it follows that
\[
ab\alpha\beta
<
(ab)^{1-1/p}[m(m+1)]^{1/p}d^2.
\]
By \eqref{eq:ab-bound} and \(p>1\),
\[
(ab)^{1-1/p}[m(m+1)]^{1/p}
<
m(m+1).
\]
Consequently, $ab\alpha\beta<m(m+1)d^2$, contradicting \eqref{eq:bad-cross}.

Thus Case $2$ also leads to a contradiction.
We have proved that in a maximizing graph, no two of
\(n_1,\ldots,n_r\) can differ by at least \(2\). Therefore
\[
\lambda^{(p)}
(K_{t-1}\vee K_r(n_1,\ldots,n_r))
\le
\lambda^{(p)}(K_{t-1}\vee T_{n-t+1,r}).
\]
Equality holds if and only if $|n_i-n_j|\le 1$ for all $i,j\in[r]$.
\end{proof}

\end{document}